\documentclass[a4paper,12pt]{article}
\author{ M. Ciavarella, L. Terracini}

\usepackage{amssymb}

\newtheorem{proposition}{Proposition}[section]

\newtheorem{corollary}{Corollary}[section]
\newtheorem{lemma}{Lemma}[section]

\newcommand{\s}{\times}
\newcommand{\QQ}{{\bf Q}}

\newcommand{\A}{{\bf A}}
\newcommand{\RR}{{\bf R}}

\newcommand{\XX}{{\bf X}}
\newcommand{\NN}{{\bf N}}

\newcommand{\OO}{\mathcal O}

\newcommand{\mA}{\mathcal A}
\newcommand\qed{\hfill \rule{6pt}{6pt}}
\newenvironment{proof}{\noindent{\bf Proof\ }\hspace{1pt}}{\hspace{-5pt}\qed\par\bigskip}

\newcommand{\ZZ}{{\bf Z}}
\def\modulo{{\rm mod}}

\def\det{{\rm det}}

\def\tr{{\rm tr}}
\def\det{{\rm det}}

\begin{document}
{}
\title{Some explicit constructions of integral structures in quaternion algebras}
\maketitle

\begin{abstract} Let $B$ be an undefined  quaternion algebra over $\QQ$. Following the explicit chacterization of some Eichler orders in $B$ given by Hashimoto,  we define explicit embeddings of these orders in some local rings of matrices; we describe the two natural inclusions of an Eichler order of leven $Nq$ in an Eichler order of level $N$. Moreover we provide a basis for a chain of Eichler orders in $B$ and prove results about their intersection. 
\end{abstract}
\section{Introduction}\label{intro}
The aim of this work is to give an explicit description of the quaternion algebras over $\QQ$ and of some of their Eichler orders. Let $B$ be a quaternion algebra over $\QQ$ of discriminant $\Delta$ and let $B_q=B\otimes_{\QQ}\QQ_q$ be its  localization at the prime number $q$. It is well known that if $q$ is a unramified place, then there is an isomorphism between $B_q$ and ${\rm M}_2(\QQ_q)$; if $B$ is ramified at $q$  then $B_q$ can be represented as a subalgebra of ${\rm M}_2(\QQ_{q^2} )$, where $\QQ_{q^2}$ denotes the quadratic unramified extension of $\QQ_q$, as described in \cite{Pizer76}.   In the general literature on quaternion algebras Eichler orders are defined by using these local isomorphisms. In \cite{Hashimoto95} an explicit definition of an Eichler orders $R(N)$ of level $N$ is given.  The author fixes a representation of the quaternion algebra $B$ as a pair $\{-\Delta N,p\}$ and gives a basis of the Eichler order $R(N)$ depending on this representation. This construction provides a very useful tool for working with Eichler orders. However, for our purposes, it has the limitation of not respecting the natural inclusion of an Eichler order of level $M$ in an Eichler order of level $N$ for $N$ dividing $M$. Starting from the work of Hashimoto, we then provide an explicit description of Eichler orders $R(N)$ and $R(Nq)$, and of the two natural inclusion maps $R(Nq)\to R(N)$.\par
More precisely, for any prime number $q$, we will describe  an isomorphisms  $\varphi_q$ between  $B_q$  and the corresponding matrix algebra and we will write the image of $R_q(N)$ under $\varphi_q$. We characterize two copies of $R(Nq)$ in $R(N)$ by using these local isomorphisms, and we define a basis for each of them in terms of a basis of $R(N)$. \\
As in \cite{Terracini03} we will consider the quaternionic analogue of the congruence groups $\Phi(N)$; we will express them by using  our characterization of Eichler ordes and we will prove some initial results for these groups.

Our interest in Eichler orders and groups $\Phi(N)$, arises from a difficulty encountered in some previous work on Galois representations and Hecke algebras arising from quaternionic groups \cite{Terracini03}, \cite{miri2006}: an analogue for Shimura curves of  Ihara's lemma (which holds for modular curves) is missing. We briefly give a sketch of this open problem; for a deep overview of the status of art  see \cite {CiavarellaTerracini}. \par

For any integer number $N$, $\Phi(N)$ is defined as $(GL_2^+(\RR)\times( R(N)\otimes\hat\ZZ)^\times )\cap B^\times$.  Let we consider the Shimura curves $\XX(N)$ and $\XX(Nq)$ coming from $\Phi(N)$ and $\Phi(Nq)$ respectively, where $q$ is a prime number such that $q\not|\Delta$. There are two injective maps from $\Phi(Nq)$ in $\Phi(N)$: the natural inclusion and  the  coniugation by a certain element $\delta_q\in B^\s$. These maps naturally induce degeneracy maps on cohomology; their direct sum provides a map $\alpha:H^1(\XX(N))^2\to H^1(\XX(Nq))$ where cohomology has coefficients in the ring of integers of a suitable finite extension of $\QQ_\ell$ for a fixed prime $\ell$. The conjecture in  \cite{CiavarellaTerracini} asserts that $\alpha$ is injective with cokernel torsion free.   

\section{Preliminaries and notations}\label{sect:preliminari}

Let $B$ be an indefinite quaternion algebra over $\QQ$ of dicriminant $\Delta=p_1...p_t$ with $t$ a even number. We will denote by $\left(\frac{*}{*}\right)$ the Legendre symbol and by $(* ,* )_q$ the Hilbert symbol at $q$ \cite{Serre70}. Let $N$ be a positive integer prime to $\Delta$ and $p$ be a prime number such that:
\begin{itemize}
\item $p\ \equiv\ 1\ \modulo\ 4$ and $p\equiv \left\{ \begin{array}{ll}
5\ \modulo\ 8 & \textrm{ if $2|\Delta$}\\
1\ \modulo\ 8 & \textrm{ if $2|N$}
\end{array} \right.$;
\item $\left(\frac{p}{p_i}\right)=-1$ for each $p_i\not=2$;
\item $\left(\frac{p}{q}\right)=1$ for each odd prime factor $q$ of $N$.
\end{itemize} 
We observe that the last condition  implies that $p$ is a square in $\ZZ_q$ for any $q$ prime factor of $N$; since $p$ is not a square in $\ZZ_p$, then $p$ does not divide $N$.\\
\noindent Hashimoto \cite{Hashimoto95} shows that then $B\simeq\{-\Delta N,p\}$ (with the notations of \cite{Vigneras80}). This means that $B$ can be expressed as $B(N,p)=\QQ+\QQ i+\QQ j+\QQ k$ where $i^2=-\Delta N$, $j^2=p$, $k=ij=-ji$.  Moreover by Theorem 2.2 of \cite{Hashimoto95}, an Eichler order  of level $N$ of $B$ can be expressed as the $\ZZ$-lattice $R(N)=\ZZ e_1+\ZZ e_2+\ZZ e_3+\ZZ e_4$ with  $$e_1=1,\  e_2=\frac{1+j}{2},\ e_3=\frac{i+k}{2},\ e_4=\frac{a\Delta Nj+k}{p}$$ where $a\in\ZZ$ satisfies $a^2\Delta N+1\equiv\ 0\ \modulo\ p$. \\
We observe that $i,j,k$ depend on the choice of $N$ and $p$;  in the sequel, 
whenever will be necessary to express the dependece on $N$ we will write $i^N,j^N,k^N$ instead of $i,j,k$ and $e_1^N,e_2^N,e_3^N,e_4^N$ instead of $e_1,e_2,e_3,e_4$. 


\noindent We consider $R=R(1)$; then $R$ is a maximal order in $B$. For any prime number $q$, let we denote $B_q=B\otimes_\QQ\QQ_q$ and $R_q(N)=R(N)\otimes_\ZZ\ZZ_q$.  \\
We start with a simple lemma which will be useful in the sequel.
\begin{lemma}\label{isom}
Let $K$ be a field and let $B_1,\ B_2$ be two quaternion algebras over $K$. If there exist a non-zero homomorphism $\varphi:B_1\to B_2$ then $\varphi$ is an isomorphism.
\end{lemma}
 \proof\\
 Since $B_1$ is a central simple algebra, it does not have non-trivial bilateral ideals  so that $\varphi$ is injective, Then the dimension $dim_K(\varphi(B_1))=4$ and $\varphi$ is an isomorphism.
 \qed
 \begin{corollary}\label{isomcor}
 Let $K$ be a field and $B_1,\ B_2$ be two quaternion algebras over $K$. We represent $B_1$ as $B_1=K+Ki+Kj+Kk$ with $i^2,j^2\in K$ and  $k=ij=-ji$.  Let $\varphi:B_1\to B_2$ be a $K$-linear map such that $$\varphi(1)=1,\ \ \varphi(i)^2=i^2,\ \ \varphi(j)^2=j^2,\ \ \varphi(k)=\varphi(i)\varphi(j)=-\varphi(j)\varphi(i).$$ Then $\varphi$ is an isomorphism of $K$-algebras. 
 \end{corollary}
We will work with $K=\QQ$ or $K=\QQ_q$ for any place $q$ including $\infty$. We observe that to define in an explicit way an isomorphism of $K$-algebras $\varphi:B_q^N\to B'$ it is enough to define the values $\varphi(i),\ \varphi(j)$ such that $\varphi(i)^2=-\Delta N$, $\varphi(j)^2=p$ and $\varphi(i)\varphi(j)=-\varphi(j)\varphi(i)$. If we put $\varphi(1)=1$, $\varphi(k)=\varphi(i)\varphi(j)$ and if we extend the map by $K$-linearity, then by Corollary \ref{isomcor},  $\varphi$ is a well defined isomorphism of $K$-algebras.

\section{The case of $ M_2(\QQ)$} \label{D=1}
If $\Delta=1$  then $B$ can be represented as $B(N,1)=\{-N,1\}$ where $N$ is any positive integer.   It is well known that there is an isomorphism $\varphi^N:B\to M_2(\QQ)$ such that the image of the maximal order $R$ is $M_2(\ZZ)$. Let we explicitly describe such an isomorphism. We consider the $\QQ$-linear map $\varphi^N$ defined as follows:
$$\varphi^N(i)=\left(
\begin{array}
[c]{cc}%
0 & -1\\
N & 0
\end{array}
\right)\ \ \textrm{and}\ \ \ \varphi^N(j)=\left(
\begin{array}
[c]{cc}%
-1 & 0\\
0 & 1
\end{array}
\right).$$ Then $$\varphi^N(i)^2=-NI\ \ \ \ \ \varphi^N(j)^2=I$$ where $I$ is the identity $2\s 2$ matrix, 
\begin{eqnarray}
\varphi^N(k)=\varphi^N(i)\varphi^N(j)&=&\left(
\begin{array}
[c]{cc}%
0 & -1\\
N & 0
\end{array}
\right)\left(
\begin{array}
[c]{cc}%
-1 & 0\\
0 & 1
\end{array}
\right)=\left(
\begin{array}
[c]{cc}%
0 & -1\\
-N & 1
\end{array}
\right)\nonumber\\
&=&-\varphi^N(j)\varphi^N(i).\nonumber
\end{eqnarray}
It results that for  any element  $x+yi+zj+tk\in B(N,1)$ with $x,y,z,t\in\QQ$ $$\varphi^N(x+iy+jz+kt)=\left(
\begin{array}
[c]{cc}%
x-z & -y-t\\
N(y-t) & x+z
\end{array}
\right)$$ and by Corollary \ref{isomcor} the map $\varphi^N:B(N,p)\to M_2(\QQ)$ is an isomorphism.
The image of the basis of the Eichler order $R(N)$ is:
$$\varphi^N(e_1)=\left(
\begin{array}
[c]{cc}%
1 & 0\\
0 & 1
\end{array}
\right)$$
$$\varphi^N(e_2)=\varphi^N\left(\frac{1+j}{2}\right)=\frac{1}{2}\left[\left(
\begin{array}
[c]{cc}%
1 & 0\\
0 & 1
\end{array}
\right)+\left(
\begin{array}
[c]{cc}%
-1 & 0\\
0 & 1
\end{array}
\right)\right]=\left(
\begin{array}
[c]{cc}%
0 & 0\\
0 & 1
\end{array}
\right)$$

$$\varphi^N(e_3)=\varphi^N\left(\frac{i+k}{2}\right)=\frac{1}{2}\left[\left(
\begin{array}
[c]{cc}%
0 & -1\\
N & 0
\end{array}
\right)+\left(
\begin{array}
[c]{cc}%
0 & -1\\
-N & 0
\end{array}
\right)\right]=\left(
\begin{array}
[c]{cc}%
0 & -1\\
0 & 0
\end{array}
\right)$$

$$\varphi^N(e_4)=\varphi^N(Nj+k)=\left(
\begin{array}
[c]{cc}%
-N & 0\\
0 & N
\end{array}
\right)+\left(
\begin{array}
[c]{cc}%
0 & -1\\
-N & 0
\end{array}
\right)=\left(
\begin{array}
[c]{cc}%
-N & -1\\
-N & N
\end{array}
\right)$$ and for any element $xe_1+ye_2+ze_3+te_4\in R(N)$ with $x,y,z,t\in\ZZ$ $$\varphi^N(xe_1+ye_2+ze_3+te_4)= \left(
\begin{array}
[c]{cc}%
x-Nt & -z-t\\
-Nt & x+y+Nt
\end{array}
\right).$$
We observe that if $N>1$, the reduced discriminant $\sqrt{|\det(\tr(\varphi^N(e_k)\varphi^N(e_h)))|}$ for $h,k=1,...,4$ of $\varphi^N(R(N))$ is $N$ so that   the image of $R(N)$ via $\varphi^N$ is $$\varphi^N(R(N))=\left\{\gamma\in M_2(\ZZ)\ |\ \gamma\equiv\left(
\begin{array}
[c]{cc}%
* & *\\
0 & *
\end{array}
\right)\ \modulo\ N\right\}.$$ 
If $N=1$ then $R(1)$ is a maximal order of $B=B^1$, the reduced discriminant of $\varphi^1(R)$ is: 
\begin{equation}
\sqrt{\left|\det\left(\begin{array}{cccc}
2 & 1 & 0 & 0\\
1 & 1 & 0 & 1\\
0 & 0 & 0 & 1\\
0 & 1 & 1 & 4
\end{array}\right)\right|}=1
\end{equation} and its image via $\varphi^1$ is $M_2(\ZZ)$.

\subsection{An isomorphism between  $B(N,1)$ and $B(M,1)$}
Let $B$ be a quaternion algebra of discriminant $1$ and let $B(N,1)=\QQ+\QQ i^N+\QQ j^N+\QQ k^N$ and $B(M,1)=\QQ+\QQ i^M+\QQ j^M+\QQ k^M$ be two representations of $B$ where $N$ and $M$ are as in Section \ref{sect:preliminari}.  We will write an isomorphism $\Psi_N^M:B(N,1)\to B(M,1)$. 

\noindent We define $\Psi_N^{M}$ as the composite $(\varphi^M)^{-1}\circ \varphi^N$:
 $$\Psi_N^{M}(i^N)=i^M\left(\frac{M+N}{2M}\right)+k^M\left(\frac{M-N}{2M}\right)$$
$$\Psi_N^{M}(j^N)=j^M$$
$$\Psi_N^M(k^N)=i^M\left(\frac{M-N}{2M}\right)+k^M\left(\frac{M+N}{2M}\right)$$
\begin{proposition}
If $M$ is an integer such that $M|N$, then $\Psi_N^M(R(N))\subset R(M)$.
\end{proposition}
\proof
Let $N=SM$ with $S\in\NN$. Then $\Psi_N^M(e_1^N)=e_1^M$, $\Psi_N^M(e_2^N)=e_2^M$, $\Psi_N^M(e_3^N)=e_3^M$ and $\Psi_N^M(e_4^N)=(1-S)e_3^M+Se_4^M$.
\qed

\section{The case of discriminant $>1$}\label{maps}
We fix a prime $p$ and a positive integer $N$ as in Section \ref{sect:preliminari}. We represent the quaternion algebra  $B$ of discriminant $\Delta$ as  $B(N,p)=\{-\Delta N,p\}=\QQ+\QQ i+\QQ j+\QQ k$. For each prime $q$ we want to identify $B_q$ to a ring of matrices, in such a way that the integer structure is preserved. Let we denote by $\mathcal R_q(N)$ the subring of $M_2(\ZZ_q)$ containing all the matrices of the form $\left\{\left(
\begin{array}
[c]{cc}%
\ZZ_q & \ZZ_q\\
N\ZZ_q & \ZZ_q
\end{array}
\right)\right\}$. We observe that if $q\not|N$ then $\mathcal R_q(N)=M_2(\ZZ_q)$. \noindent We recall that every local Eichler orders  of level $N$ in $M_2(\QQ_q)$ is isomorphic to $\mathcal R_q(N)$ and its reduced discriminant is equal to $\Delta N$. \\

\noindent We will deal separately with the cases of unramified places and of ramified places.\\

\subsection{The isomorphism at the non-Archimedean unramified places}

In this section let $q$ be a prime number  such that $q\not|\Delta$; since  at $q$ the quaternion algebra $B(N,p)$ is not ramified, the Hilbert symbol is 
\begin{equation}\label{Hilbert}
1=(-\Delta N,p)_q.
\end{equation}
 We shall define an isomorphism $\varphi^{(N,p)}_q:B(N,p)_q\widetilde{\to} M_2(\QQ_q)$ such that $\varphi_q^{(N,p)}(R_q(N))=\mathcal R_q(N).$ To make easier the notation we will write $\varphi_q^N$ instead of $\varphi_q^{(N,p)}$.

\subsubsection{The isomorphism at places $q$ not  dividing $\Delta p$ such that $p$ is not a square in $\ZZ_q^\s$}\label{primi|N}

We consider the case $q\not |\Delta p$ such that $\left(\frac{p}{q}\right)=-1$ (we observe that the last condition excludes the cases $q=2$ and $q|N$).
This hypotheses on $q$ assure that $p$ is not a square in $\ZZ_q^\s$, thus $\QQ_q(\sqrt{p})$ is a quadratic extension of $\QQ_q$ and by the identity (\ref{Hilbert}), the prime $-\Delta N$ is the norm of a unit of $\QQ_q(\sqrt{p})$. We write $-\Delta N=x^2-py^2$ with $x,y\in\ZZ_q$ and we define  $\varphi_q^N$ as follows:

$$\varphi_q^N(i)=\left(
\begin{array}
[c]{cc}%
x & -py\\
y & -x
\end{array}
\right)\ \ \ \ \varphi_q^N(j)=\left(
\begin{array}
[c]{cc}%
0 & p\\
1 & 0
\end{array}
\right).$$
Then $$\varphi_q^N(i)^2=-\Delta NI\ \ \ \ \ \varphi_q^N(j)^2=pI$$ $$\varphi_q^N(k)=\varphi_q^N(i)\varphi_q^N(j)=\left(
\begin{array}
[c]{cc}%
-py & px\\
-x & py
\end{array}
\right) =-\varphi_q^N(j)\varphi_q^N(i).$$
It results that for any $h=\alpha+\beta i+\gamma j+\delta k\in B_q(N,p)$
$$\varphi_q^N(h)=\left(
\begin{array}
[c]{cc}%
\alpha+\beta x-\delta py & -\beta py+\gamma p+\delta xp\\
\beta y+\gamma-\delta x & \alpha-\beta x+\delta py
\end{array}
\right) $$ and by Corollary \ref{isomcor}, $\varphi_q^N:B_q\to M_2(\QQ_q)$ is an isomorphism.\\
The image of the basis of the local Eichler order is:
$$\varphi_q^N(e_1)=\left(
\begin{array}
[c]{cc}%
1 & 0\\
0 & 1
\end{array}
\right)$$
$$\varphi_q^N(e_2)=\frac{1}{2}\left(
\begin{array}
[c]{cc}%
1 & p\\
1 & 1
\end{array}
\right)\in M_2(\ZZ_q)\ \textnormal{since}\ q\not=2$$
$$\varphi_q^N(e_3)=\frac{1}{2}\left(
\begin{array}
[c]{cc}%
x-py & p(x-y)\\
y-x & -x+py
\end{array}
\right)\in M_2(\ZZ_q)\ \textnormal{since}\ q\not=2$$
$$\varphi_q^N(e_4)=\left(
\begin{array}
[c]{cc}%
-y & a\Delta N+x\\
\frac{a\Delta N-x}{p} & y
\end{array}
\right)\in M_2(\ZZ_q)\ \textnormal{since}\ q\not=p.$$
The reduced discriminant of $\varphi_q^N(R_q(N))$ is $\Delta N$. So $\varphi_q^N(R_q(N))=M_2(\ZZ_q)$. For any element $g=\alpha e_1+\beta e_2+\gamma e_3+\delta e_4\in R_q(N)$
\begin{equation}\label{q^2|N}
\varphi_q^N(g)=\left(
\begin{array}
[c]{cc}%
\alpha+\frac{\beta}{2}+\frac{\gamma}{2}(x-py)-\delta y & \beta\frac{p}{2}+\frac{\gamma p}{2}(x-y)+\delta(a\Delta N+x)\\
\frac{\beta}{2}+\frac{\gamma}{2}(y-x)+\delta\frac{a\Delta N-x}{p} & \alpha+\frac{\beta}{2}+\frac{\gamma}{2}(py-x)+\delta y
\end{array}
\right)
\end{equation}

\subsubsection{The isomorphism at primes $q$ such that $p$ is a square in $\ZZ_q^\s$}\label{q=2}
We consider the primes $q\not|\Delta$ such that $\left(\frac{p}{q}\right)=1$. We observe that this hypothesis excludes the case $p=q$ and includes $q|N$ and $q=2$ (in fact if $q=2$ then  by hypothesis $p\equiv 1\ \modulo\ 8$ and by (\cite{Serre70}, II, \S3) $p$ is a square in $\ZZ_2^\s$).
We define the $\QQ_q$-linear map $\varphi_q^N$ as follows:
$$\varphi_q^N(i)=\left(
\begin{array}
[c]{cc}%
0 & 1\\
-\Delta N & 0
\end{array}
\right)\ \ \textnormal{and}\ \ \varphi_q^N(j)=\left(
\begin{array}
[c]{cc}%
-\sqrt p & 0\\
0 & \sqrt p
\end{array}
\right)$$ where $\sqrt p$ is an element $\omega$ in $\ZZ_q^\s$ such that $\omega^2=p$. Then
$$\varphi_q^N(i)^2=-\Delta N I\ \ \ \ \ \varphi_q^N(j)^2=pI$$ 
$$\varphi_q^N(k)=\varphi_q^N(i)\varphi_q^N(j)=\left(
\begin{array}
[c]{cc}%
0 & \sqrt p\\
\Delta N\sqrt p & 0
\end{array}
\right)=-\varphi_q^N(j)\varphi_q^N(i).$$ 
It results that for any element $\alpha+\beta i+\gamma j+\delta k\in B_q$ 
$$\varphi_q^N(\alpha+\beta i+\gamma j+\delta k)= \left(
\begin{array}
[c]{cc}%
\alpha-\gamma\sqrt p & \beta+\delta\sqrt p \\
\Delta N(-\beta+\delta\sqrt p) & \alpha+\gamma\sqrt p
\end{array}
\right)$$ 
and by Corollary \ref{isomcor}, $\varphi_q^N:B_q(N,p)\to M_2(\QQ_q)$ is an isomorphism.\\
The image of a basis of the local Eichler order $R_q(N)$ is:
$$\varphi_q^N(e_1)=I$$
$$\varphi_q^N(e_2)=\left(
\begin{array}
[c]{cc}%
\frac{1-\sqrt p}{2} & 0\\
0 & \frac{1+\sqrt p}{2}
\end{array}
\right)\in \mathcal R_q(N)$$
$$\varphi_q^N(e_3)=\left(
\begin{array}
[c]{cc}%
0 & \frac{1+\sqrt p}{2}\\
\frac{\Delta N(\sqrt p-1)}{2} & 0
\end{array}
\right)\in \mathcal R_q(N)$$
$$\varphi_q^N(e_4)=\frac{1}{\sqrt p}\left(
\begin{array}
[c]{cc}%
-a\Delta N & 1\\
\Delta N & a\Delta N
\end{array}
\right)\in \mathcal R_q(N)$$
The reduced discriminant of $\varphi_q^N(R_q(N))$ is $\Delta N$, so $\varphi_q^N(R_q(N))=\mathcal R_q(N)$. 
Then, for any element $g=\alpha e_1+\beta e_2+\gamma e_3+\delta e_4\in R_q(N)$  
\begin{equation}\label{eldue}
\varphi_q^N(g)=\\
\left(
\begin{array}
[c]{cc}%
\alpha+\frac{(1-\sqrt p)\beta}{2}-\frac{\delta a\Delta N}{\sqrt p} & \gamma\frac{(1+\sqrt p)}{2}+\frac{\delta}{\sqrt p} \\
\gamma\Delta N\frac{\sqrt p-1}{2}+\frac{\delta\Delta N}{\sqrt p} & \alpha+\frac{(1+\sqrt p)\beta}{2}+\frac{\delta a\Delta N}{\sqrt p}
\end{array}
\right).
\end{equation}
We observe that in this case we accept that $2|N$.

\subsubsection{The isomorphism at $p$}\label{isoq=p}

\noindent  If $q=p$ then $1=(-\Delta N,p)_p=\left(\frac{-\Delta N}{p}\right)$ (\cite{Serre70}, II, \S3). So $-\Delta N$ is a square in $\ZZ_p^\s$. We recall that $a\in\ZZ$ was choosen in Section \ref{sect:preliminari} in such a way that $a^2\Delta N+1\equiv 0\ \modulo\ p$. Let we denote by $\sqrt{-\Delta N}$ the square root of $-\Delta N$ in $\ZZ_p^\s$ such that $a\sqrt{-\Delta N}\equiv\ -1\ \modulo\ p$. Then  the following identity holds:
\begin{equation}\label{cN}
(a\Delta N-\sqrt{-\Delta N})=\sqrt{-\Delta N}(a\sqrt{-\Delta N}-1)\equiv 0\ \modulo\ p.
\end{equation}

\noindent We define the $\QQ_p$-linear map $\varphi_p^N$ as follows:  
$$\varphi_p^N(i)=\left(
\begin{array}
[c]{cc}%
-\sqrt{-\Delta N} & 0\\
0 & \sqrt{-\Delta N}
\end{array}
\right)\ \ \textnormal{and}\ \ \varphi_p^N(j)=\left(
\begin{array}
[c]{cc}%
0 & 1\\
p & 0
\end{array}
\right).$$
Then
$$\varphi_p^N(i)^2=-\Delta NI\ \ \ \ \ \varphi_p^N(j)^2=pI$$ 
$$\varphi_p^N(k)=\varphi_p^N(i)\varphi_p^N(j)=\left(
\begin{array}
[c]{cc}%
0 & -\sqrt{-\Delta N}\\
p\sqrt{-\Delta N} & 0
\end{array}
\right)=-\varphi_p^N(j)\varphi_p^N(i).$$
It results that  for any element $\alpha+\beta i+\gamma j+\delta k\in B_p(N,p)$ $$\varphi_p^N(\alpha+\beta i+\gamma j+\delta k)= \left(
\begin{array}
[c]{cc}%
\alpha-\beta\sqrt{-\Delta N} & \gamma-\delta\sqrt{-\Delta N}\\
\gamma p+\delta p\sqrt{-\Delta N} & \alpha+\beta\sqrt{-\Delta N}
\end{array}
\right)$$ and by Corollary \ref{isomcor}, $\varphi_p^N:B(N,p)_p\to M_2(\QQ_p)$ is an isomorphism.  It remains to show that integer structures are preserved.

\noindent The image of the basis of the local Eichler order is:
$$\varphi_p^N(e_1)=I$$
$$\varphi_p^N(e_2)=\frac{1}{2}\left(
\begin{array}
[c]{cc}%
1 & 1\\
p & 1
\end{array}
\right)\in M_2(\ZZ_p)$$
$$\varphi_p^N(e_3)=\frac{1}{2}\left(
\begin{array}
[c]{cc}%
-\sqrt{-\Delta N} & -\sqrt{-\Delta N}\\
p\sqrt{-\Delta N} & \sqrt{-\Delta N}
\end{array}
\right)\in M_2(\ZZ_p)$$
$$\varphi_p^N(e_4)=\left(
\begin{array}
[c]{cc}%
0 & \frac{a\Delta N-\sqrt{-\Delta N}}{p}\\
a\Delta N+\sqrt{-\Delta N} & 0
\end{array}
\right)\in M_2(\ZZ_p).$$ 
 
\noindent The reduced discriminant of $\varphi_p^N(R_p(N))$ is $\Delta N$ so that $\varphi_p^N(R_p(N))=M_2(\ZZ_p)$.\\
For any element $\alpha e_1+\beta e_2+\gamma e_3+\delta e_4\in R_p(N)$ the following identity holds:
\begin{eqnarray}\label{p=q} 
&\varphi_p^N(\alpha e_1+\beta e_2+\gamma e_3+\delta e_4)=&\\
&\frac{1}{2}\left(
\begin{array}
[c]{cc}%
2\alpha+\beta-\gamma\sqrt{-\Delta N} & \beta-\gamma\sqrt{-\Delta N}+\frac{2\delta}{p}(a\Delta N-\sqrt{-\Delta N}) \\
p\beta+\gamma p\sqrt{-\Delta N}+2\delta(a\Delta N+\sqrt{-\Delta N}) & 2\alpha+\beta+\gamma\sqrt{-\Delta N}
\end{array}
\right).&\nonumber
  \end{eqnarray}

\subsection{The isomorphism at the Archimedean place}

\noindent Since $B$ is an indefinite quaternion algebra over $\QQ$,  there exists an isomorphism $B_\infty\simeq M_2(\RR)$.
We define $\varphi^N_\infty$ via $$i\mapsto\left(
\begin{array}
[c]{cc}%
0 & 1\\
-\Delta N & 0
\end{array}
\right)\ \ \ \ \ j\mapsto\left(
\begin{array}
[c]{cc}%
\sqrt p & 0\\
0 & -\sqrt p
\end{array}
\right).$$
Then $$\varphi^N_\infty(i)^2=-\Delta NI\ \ \ \ \ \varphi^N_\infty(j)^2=pI$$ $$\varphi^N_\infty(k)=\varphi^N_\infty(i)\varphi^N_\infty(j)=\left(
\begin{array}
[c]{cc}%
0 & -\sqrt p\\
-\sqrt p\Delta N & 0
\end{array}
\right)=-\varphi^N_\infty(j)\varphi^N_\infty(i)$$
and by Corollary \ref{isomcor} the map $\varphi_\infty^N:B_\infty\to M_2(\RR)$ is an isomorphism.

\subsection{The isomorphism at the ramified places}
For any prime number $q$ such that $q|\Delta$ we shall define, following \cite{Pizer76}, an isomorphism $\varphi_q^N:B_q(N,p)\widetilde\to \left\{\left(
\begin{array}
[c]{cc}%
\alpha & \beta\\
q\overline\beta & \overline\alpha
\end{array}
\right)\ |\ \alpha,\beta\in \QQ_{q^2}  \right\}$ such that $$\varphi_q^N(R_q(N))=\varphi_q^N(R_q)=\left\{\left(
\begin{array}
[c]{cc}%
\alpha & \beta\\
q\overline\beta & \overline\alpha
\end{array}
\right)\ |\ \alpha,\beta\in  \ZZ_{q^2}) \right\}:=\OO_q$$
 where $ \QQ_{q^2}$ is  the quadratic unramified extension of $\QQ_q$, $\alpha\mapsto \bar\alpha$ is its non-trivial automorphism  and  $ \ZZ_{q^2}$ is its ring of integers.

\noindent We have  $-1=(-\Delta N,p)_q=(-\frac{\Delta N}{q},p)_q(q,p)_q$; since $\Delta$ is square free $(-\frac{\Delta N}{q},p)_q=1$ and $(q,p)_q=-1$. This means in particular that $p$ is not a square in $\QQ_q$ and  $-\frac{\Delta N}{q}$ is a norm of a unit of $\QQ_q(\sqrt p)$. Thus there exist $x,y\in\ZZ_q$ such that $-\frac{\Delta N}{q}=x^2-py^2=(x-\sqrt py)(x+\sqrt py)$. We can identify $\QQ_{q^2}=\QQ_q(\sqrt p)$  and $\ZZ_{q^2}=\ZZ_q(\sqrt p)$.\\
We define  $\varphi_q^N$ as follows:
$$\varphi_q^N(i)=\left(
\begin{array}
[c]{cc}%
0 & x-\sqrt py\\
q(x+\sqrt py) & 0
\end{array}
\right)$$
$$\varphi_q^N(j)=\left(
\begin{array}
[c]{cc}%
-\sqrt p & 0\\
0 & \sqrt p
\end{array}
\right).$$
Then $$\varphi^N_q(i)^2=\Delta N I\ \ \ \ \ \varphi_q^N(j)^2=pI$$
$$\varphi_q^N(k)=\varphi_q^N(i)\varphi_q^N(j)=\left(
\begin{array}
[c]{cc}%
0 & \sqrt p(x-\sqrt py)\\
-\sqrt pq(x+\sqrt py) & 0
\end{array}
\right)=-\varphi_q^N(j)\varphi_q^N(i)$$
and for any element $\alpha+\beta i+\gamma j+\delta k$ of $B_q(N,p)$ with $\alpha,\beta,\gamma,\delta\in\QQ_q$, $$\varphi_q^N(\alpha+\beta i+\gamma j+\delta k)= \left(
\begin{array}
[c]{cc}%
\alpha-\gamma\sqrt p & (\beta+\delta\sqrt p)(x-\sqrt py)\\
q(\beta-\delta\sqrt p)(x+\sqrt py) & \alpha+\gamma\sqrt p
\end{array}
\right).$$ By Corollary \ref{isomcor}, $\varphi_q^N$ is an isomorphism.

\noindent We compute the image of the local Eichler order $R_q(N)$: 
$$\varphi_q^N(e_1)=I$$
$$\varphi_q^N(e_2)=\frac{1}{2}\left(
\begin{array}
[c]{cc}%
1-\sqrt p & 0\\
0 & 1+\sqrt p
\end{array}
\right)\in\OO_q$$
$$\varphi_q(e_3)=\frac{1}{2}\left(
\begin{array}
[c]{cc}%
0 & (x-\sqrt py)(1+\sqrt p)\\
q(x+\sqrt py)(1-\sqrt p) & 0
\end{array}
\right)\in\OO_q$$
$$\varphi_q(e_4)=\left(
\begin{array}
[c]{cc}%
\frac{-a\Delta N}{p}\sqrt p & -y+\frac{x}{p}\sqrt p\\
q\left(-y-\frac{x}{p}\sqrt p\right) & \frac{a\Delta N}{p}\sqrt p
\end{array}
\right)\in\OO_q$$
and the reduced discriminant of $\varphi^N_q(R_q(N))$ is $N\Delta$.\\
For any element $\alpha e_1+\beta e_2+\gamma e_3+\delta e_4$ of $R_q(N)$ $$\varphi_q^N(\alpha e_1+\beta e_2+\gamma e_3+\delta e_4)=$$ $$\left(
\begin{array}
[c]{cc}%
\alpha+\frac{\beta}{2}-\sqrt p\left(\frac{\beta}{2}+aN\Delta\frac{\delta}{p}\right) & (x-\sqrt py)\left[\frac{\gamma}{2}+\sqrt p\left(\frac{\gamma}{2}+\frac{\delta}{p}\right)\right] \\
q(x+\sqrt py)\left[\frac{\gamma}{2}-\sqrt p\left(\frac{\gamma}{2}+\frac{\delta}{p}\right)\right] & \alpha+\frac{\beta}{2}+\sqrt p\left(\frac{\beta}{2}+aN\Delta\frac{\delta}{p}\right)
\end{array}
\right).$$

\section{Characterization of $\Phi(N)$}\label{caratterizzazione}

Let $B$ be a quaternion algebra over $\QQ$ of discriminant $\Delta$ and let  $R(N)$ be an Eichler order of level $N$ of $B$. Then $R(N)^\s\s\hat{\ZZ}$ is a compact open subgroup of the finite adelization  $B_{\A}^{\s,\infty}$ and it is possible to associate to it a discrete subgroup $\Phi(N)$ of $SL_2(\RR)$ by
 $$\Phi(N)=(GL_2^+(\RR)\s (R(N)^\s\s\hat{\ZZ}) )\cap B^\s.$$
 It is known that $\Phi(N)$ is a co-compact congruence subgroup of $SL_2(\RR)$ \cite{Vigneras80}.

\begin{lemma}\label{uno}
If we  denote by $R(N)^{(1)}$ the group of reduced norm $1$ elements of $R(N)$, then the following identity holds: $\Phi(N)=R(N)^{(1)}.$
\end{lemma}
\proof
The inclusion $\supseteq$ is trivial since $R(N)^{(1)}\subseteq B^\s$ and $R(N)^{(1)}\subseteq GL_2^+(\RR)\s (R(N)^\s\s\hat{\ZZ})$.\\
We prove the inclusion $\subseteq$. Let $\alpha$ be an element of $\Phi(N)$; then:

\begin{itemize}
\item[a)] $\alpha\in GL_2^+(\RR)\s (R(N)^\s\s\hat{\ZZ})$
\item[b)] $\alpha\in B^\s$
\end{itemize}
If we denote by $n(\alpha)$ the reduced norm of $\alpha$, then by b), $n(\alpha)$ is a rational number, which, by a), is a $p$-adic unit for every prime $p$, and positive. Thus  $n(\alpha)=1$ and $\alpha\in R(N)$. \qed

\section{Explicit description of two conjugates to $R(Nq)$ in $B(N,p)$}
In this section we will keep the usual notation and we will represent the quaternion algebra $B$ as $B(N,p)=\{-N\Delta,p\}$.

\noindent Let $q$ be a prime number such that $q\not|\Delta$. It is well known that by definition
\begin{equation}\label{R(Nq)}
R(Nq)\simeq  R(N)\cap(\varphi_q^N)^{-1}(\mathcal R_q(qN)).
\end{equation}
We shall identify $R(Nq)$ with this subgroup of $R(N)$.
\noindent Let we consider the id\`ele $\eta_q$ in $B_\A^\s$ defined by

\begin{displaymath}
\eta_q=\left\{\begin{array}{ll}
\eta_{q,\nu}=1 & \rm{if}\ \nu\not=q\\
\eta_{q,q}=(\varphi_q^N)^{-1}\left(
\begin{array}
[c]{cc}%
q & 0\\
0 & 1
\end{array}
\right) & \rm{if}\ \nu=q
\end{array}\right. 
\end{displaymath}

By strong approximation, write $\eta_q=\delta_qg_\infty u$, with $\delta_q\in B^\s$, $g_\infty\in GL_2^+(\RR)$ and $u\in(R(Nq)\otimes_\ZZ\widehat\ZZ)^\s$.\\
We observe that $$\eta_qR_q(Nq)\eta_q^{-1}=\delta_qR_q(Nq)\delta_q^{-1}=R_q(N)\cap( \varphi_q^N)^{-1}\left(
\begin{array}
[c]{cc}%
\ZZ_q & q\ZZ_q\\
N\ZZ_q & \ZZ_q
\end{array}
\right)$$ and 
\begin{equation}\label{deltaR(Nq)}
\delta_qR(Nq)\delta_q^{-1}=R(N)\cap( \varphi_q^N)^{-1}\left(
\begin{array}
[c]{cc}%
\ZZ_q & q\ZZ_q\\
N\ZZ_q & \ZZ_q
\end{array}
\right).
\end{equation}

\noindent We will give bases for $R(Nq)$ and $\delta_qR(Nq)\delta_q^{-1}$. We observe that the following  theorems are direct applications of the construction in \cite{Samuel71} \S1.5, by considering the results in the previous sections and the image via the isomorphisms $\varphi_q^N$ of a generic element of $R_q(N)$.

\begin{proposition}\label{NNqualunque}
Let $q$ be a prime number such that $q\not|\Delta p$ and $p$ is not a square in $\ZZ_q^\s$. Let $-\Delta N=x^2-py^2$ with $x,y\in\ZZ_q$. Let $c_1,c_2,c_3$ be integers  such that $$c_1\equiv (y-x)\ \modulo \ q$$ $$c_2\equiv p^{-1}\ \modulo\ q$$ $$c_3\equiv x\ \modulo\ q.$$ 
Then  a basis of $R(Nq)$ in $R(N)$ is: $$f_1=e_1,\ \ f_2=-c_1e_2+e_3,\ \ f_3=-2c_2(a\Delta N-c_3)e_2+e_4,\ \ f_4=qe_2$$
and a basis of $\delta_qR(Nq)\delta_q^{-1}$ in $R(N)$ is: $$g_1=e_1,\ \ g_2=c_1e_2+e_3,\ \ g_3=-2c_2(a\Delta N+c_3)e_2+e_4,\ \ g_4=qe_2.$$
\end{proposition}
\proof
By the results in Section \ref{primi|N} and by the equality (\ref{R(Nq)}), we see that $f_1, f_2, f_3, f_4\in R(Nq)$ and  $$\det\left(
\begin{array}{cccc}
1 & 0 & 0 & 0\\
0 & -c_1 & -2c_2(a\Delta N-c_3) & q\\
0 & 1 & 0 & 0\\
0 & 0 & 1 & 0
\end{array}
\right)=q.$$  By the results in Section \ref{primi|N} and by the equality (\ref{deltaR(Nq)}) we see that $g_1, g_2, g_3, g_4\in \delta_qR(Nq)\delta_q^{-1}$ and $$\det\left(
\begin{array}{cccc}
1 & 0 & 0 & 0\\
0 & c_1 & -2c_2(a\Delta N+c_3) & q\\
0 & 1 & 0 & 0\\
0 & 0 & 1 & 0
\end{array}
\right)=q.$$ \qed

\begin{proposition}
Let $q\not|\Delta$ be a prime number such that $p$ is a square in $\ZZ_q^\s$. A basis of $R(qN)$ in $R(N)$ is: $$f_1=e_1,\ \ f_2=e_2,\ \ f_3=e_3-ce_4,\ \ f_4=qe_4$$ where $\ZZ\ni c\equiv\ (p-\sqrt p)2^{-1}\ \modulo\ q.$ 
A basis of $\delta_qR(qN)\delta_q^{-1}$ in $R(N)$ is: $$g_1=e_1,\ \ g_2=e_2,\ \ g_3=e_3-c'e_4,\ \ g_4=qe_4$$ where $\ZZ\ni c'\equiv\ (p+\sqrt p)2^{-1}\ \modulo\ q.$  
\end{proposition}

\proof 
By the results in section \ref{q=2} and by the equality (\ref{R(Nq)}), we observe that $f_1,f_2,f_3,f_4\in R(qN)$ and $$\det\left(
\begin{array}{cccc}
1 & 0 & 0 & 0\\
0 & 1 & 0 & 0\\
0 & 0 & 1 & 0\\
0 & 0 & -c & q
\end{array}
\right)=q;$$
by the equality (\ref{deltaR(Nq)}), it is easy to verify that $g_1,g_2,g_3,g_4\in \delta_qR(qN)\delta_q^{-1}$ and $$\det\left(
\begin{array}{cccc}
1 & 0 & 0 & 0\\
0 & 1 & 0 & 0\\
0 & 0 & 1 & 0\\
0 & 0 & -c' & q
\end{array}
\right)=q.$$\qed

\begin{proposition}\label{Np}
Let $\sqrt{-\Delta N}$ be the square root of $-\Delta N$ in $\QQ_p$ such that $a\sqrt{-\Delta N}\equiv -1\ \modulo\ p$.  A basis of $R(Np)$ in $R(N)$ is: $$f_1=e_1,\ \ f_2=-c_4e_2+e_3,\ \ f_3=-2(a\Delta N+c_4)e_2+pe_4,\ \ f_4=p(Ae_2+Be_4)$$ and a basis of $\delta_pR(Np)\delta_p^{-1}$ in $R(N)$ is: $$g_1=e_1,\ \ g_2=c_4e_2+e_3,\ \ g_3=-2\frac{a\Delta N-c_4}{p}e_2+e_4\ \ g_4=pe_2$$ where $\ZZ\ni c_4\equiv \sqrt{-\Delta N}\ \modulo\ p$,  $a\in\ZZ$ is such that $a^2\Delta N+1\equiv 0\ \modulo\ p$ and $A,B\in\ZZ$ are such that $Ap+2B(a\Delta N+c_4)=1$.\\
\end{proposition}

\proof
We first observe that if we fix $\sqrt{-\Delta N}$  the square roots in $\QQ_p$ such that $a\sqrt{-\Delta N}\equiv -1\ \modulo\ p$, then $p|(a\Delta N-c_4)$  and $p\not|(a\Delta N+c_4)$. Then the existence of $A,B\in\ZZ$  such that $Ap+2B(a\Delta N+c_4)=1$ is ensured.\\
By the results in Section \ref{isoq=p} and by the equality (\ref{R(Nq)}), we observe that $f_1,f_2,f_3,f_4\in R(Np)$ and $$\det\left(
\begin{array}{cccc}
1 & 0 & 0 & 0\\
0 & -c_4 & -2(a\Delta N+c_4) & pA\\
0 & 1 & 0 & 0\\
0 & 0 & p & pB
\end{array}
\right)=p;$$
by the equality (\ref{deltaR(Nq)}), we see that 
$g_1,g_2,_3,g_4\in\delta_pR(Np)\delta_p^{-1}$ and $$\det\left(
\begin{array}{cccc}
1 & 0 & 0 & 0\\
0 & c_4 & -2\frac{a\Delta N-c_4}{p} & p\\
0 & 1 & 0 & 0\\
0 & 0 & 1 & 0
\end{array}
\right)=p.$$\qed

\section{Description of the isomorphism $\Psi_N^{M}$}\label{isomo}

Let we fix $\Delta$ as in Section \ref{sect:preliminari};  by the classification theorem, up to isomorphism there exists only one quaternion algebra $B$ over $\QQ$ with discriminant  $\Delta$. Let $B(M,p)=\{-\Delta M,p\}$
and $B(N,p)=\{-\Delta N,p\}$
be two representations of $B$, with $p,N,M$ as in Section \ref{sect:preliminari}. Then there is an isomorphism $\Psi_N^M:B(N,p)\to B(M,p)$. This implies that there exists an element $h\in B(M,p)$ such that $h^2=-N\Delta$; we observe that $h$ is of the form  $h=i^M\beta+j^M\gamma+k^M\delta$ where $(\beta, \gamma, \delta)\in\QQ^3$ is a solution of the equation $$M\Delta \beta^2-p\gamma^2-p\Delta M\delta^2=N\Delta.$$
\begin{lemma}\label{quadratic}
Let $f$ be the quadratic form on $\QQ$ defined as $f=M\beta^2-pM\delta^2$; then $f$ represents $N$.
\end{lemma}
\proof
By the Hasse-Minkowski theorem (see for example \cite{Serre70}), $f$ represents $N$ in $\QQ$ if and only if $f$ represents $N$ in $\QQ_\ell$ at any place $\ell$, that is $(N,p)_\ell=(M,p)_\ell$ for any prime number $\ell$.\\ 
We write: $N=\ell^au$, $p=\ell^bv$, $\epsilon(\ell)\equiv\frac{\ell-1}{2}\ \modulo\ 2$.\\
If $\ell\not=2$ then $$(N,p)_\ell=(-1)^{ab\epsilon(\ell)}\left(\frac{u}{\ell}\right)^b\left(\frac{v}{\ell}\right)^a$$
\begin{itemize}
\item If $\ell\not|pN$ then $(N,p)_\ell=1$.
\item If $\ell=p$ then $\epsilon(p)=0,\ v=1,\ b=1$ so $(N,p)_p=\left(\frac{u}{p}\right)^b\left(\frac{v}{p}\right)^a=\left(\frac{u}{p}\right)$. By the hypothesis on the prime factors $q$ of $N$, by the law  of  recipocity and since $p\equiv 1\ \modulo\ 4$, $$\left(\frac{u}{p}\right)=\prod_{q|u}\left(\frac{q}{p}\right)=\prod_{q|u}\left(\frac{p}{q}\right)(-1)^{(q-1)(p-1)/4}=1.$$
\item If $\ell|N$ and $\ell\not=p$ then $b=0, v=p$ so $$(N,p)_\ell=\left(\frac{p}{\ell}\right)^a=1$$ by the hypothesis on the prime factors of $N$.
\end{itemize}
If $\ell=2$ then $b=0$ and $v=p$; we know that $$(N,p)_2=(-1)^{\epsilon(u)\epsilon(v)+a\omega(p)+b\omega(u)}$$ where $\epsilon(v)=0, \omega(p)\equiv\frac{p^2-1}{8}\ \modulo\ 2$. So $$(N,p)_2=(-1)^{a\omega(p)}.$$
\begin{itemize}
\item if $a=0$ then $(N,p)_2=1$;
\item if $a\not=0$ then $2|N$ and $p\equiv 1\ \modulo\ 8$. So $\omega(p)=0$ and $(N,p)_2=1.$
\end{itemize}
{}

\noindent Since $N$ and $M$ satisfy the same hypotheses, then $(N,p)_\ell=(M,p)_\ell=1$ for any prime number $\ell$.
\qed

\noindent We define the $\QQ$-linear map $\Psi_N^M:B(N,p)\to B(M,p)$ as:  $$\Psi_N^M(i^N)=h,\ \ \Psi_N^M(j^N)=j^M$$ where $h=\beta i^M+\delta k^M$ with $(\beta, \gamma)\in\QQ^2$  solution of $M\beta^2-pM\delta^2=N$ (by Lemma \ref{quadratic} such an element exists).   Then $\Psi_N^M(i^N)^2=-N\Delta$, $\Psi_N^M(j^N)^2=p$ and $$\Psi_N^M(i^N)\Psi_N^M(j^N)=hj^M=k^M\beta+i^Mp\delta=-\Psi_N^M(j^N)\Psi_N^M(i^N).$$  By the Corollary \ref{isomcor}, the map $\Psi_N^M$ is an isomorphism. \\

\noindent We observe that if $N=MS$  then 
\begin{equation}\label{equazioneS}
\beta^2-p\delta^2=S 
\end{equation}
so $S$ is the norm  of an element $\beta+\sqrt p\delta$ of the ring of integer 
 $$\OO=\left\{\frac{1}{2}(a+\sqrt pb)\ :\ a,b\in\ZZ\ \textnormal{with the same parity} \right\}$$ of $\QQ(\sqrt p)$.

\noindent We denote by $a_M, a_N$ the integer numbers as in Section \ref{sect:preliminari}, such that $a_M^2\Delta M+1\equiv\ 0\ \modulo\ p$ and  $a_N^2\Delta N+1\equiv\ 0\ \modulo\ p$.
\begin{lemma}\label{a_M}
If $N=MS\in\NN$ then we can choose $\beta\in\ZZ\left[\frac{1}{2}\right]$ satisfying the identity (\ref{equazioneS}) such that $a_M\equiv\ a_N\beta\ \modulo\ p$. 
\end{lemma}
\proof
By definition of $a_N, a_M$, since $p\not|\Delta M$, we find that $a_N^2S-a_M^2\equiv 0\ \modulo\ p$, that is by (\ref{equazioneS}) $a_N^2\beta^2-a_M^2\equiv 0\ \modulo\ p$. In particular $a_N\beta-a_M\equiv 0\ \modulo\ p$ or $a_N\beta+a_M\equiv 0\ \modulo\ p$.  If we are in the second situation, then we can take $-\beta$ instead of $\beta$, so we have that $a_M\equiv a_N\beta\ \modulo\ p$.
\qed

\noindent In the sequel when $M|N$ we choose $a_M$ as in the above lemma.

\begin{proposition}
Let $B(N,p)$ and $B(M,p)$ be two representations of the quaternion algebra $B$ defined over $\QQ$ with discrminant $\Delta$. Let we consider the isomorphism $\Psi_N^M:B(N,p)\to B(M,p)$ defined above. If $M|N$ then $\Psi_N^M(R(N))\subset R(M)$. 
\end{proposition}
\proof
Let $N=MS$ where $S\in\NN$.
We recall that $p\equiv\ 1\ \modulo\ 4$ and we verify that $\Psi_N^M(e_\ell^N)\in R(M)$ for $\ell=1,2,3,4$.\\
By definition of $\Psi_N^M$:
$$\Psi_N^M(e_1^N)=1=e_1^M$$
$$\Psi_N^M(e_2^N)=\frac{1+j^M}{2}=e_2^M$$
$$\Psi_N^M(e_3^N)=A_3e_1^M+B_3e_2^M+C_3e_3^M+D_3e_4^M$$ where
$A_3=\frac{1}{2}\delta(1-p)a_M\Delta M\in\ZZ$, $B_3=\delta(p-1)a_M\Delta M\in\ZZ$, $C_3=\delta p+\beta\in\ZZ$ and $D_3=\delta p\frac{1-p}{2}\in\ZZ$.
$$\Psi_N^M(e_4^N)=A_4e_1^M+B_4e_2^M+C_4e_3^M+D_4e_4^M$$ where
$B_4=-2A_4=\frac{2}{p}\left[\Delta M(a_NS-a_M\beta+p\delta a_M\right]$, $C_4=2\delta\in\ZZ$ and $D_4=\beta-p\delta\in\ZZ$. We observe that $B_4\in\ZZ$ (and $A_4\in\ZZ$), infact by Lemma \ref{a_M}:
\begin{eqnarray}
a_NS-a_M\beta&\equiv\ &a_NS-a_N\beta^2\ \modulo\ p\nonumber\\
&\equiv\ &a_NS-a_N(S+p\delta^2)\ \modulo\ p\nonumber\\
&\equiv\ &0\ \modulo\ p\nonumber
\end{eqnarray}
\qed

\section{Some properties of the Eichler orders}

By using the local isomorphisms given in Section  \ref{maps}, we will prove some new results for the Eichler orders.  Let $B$ be a quaternion algebra over $\QQ$ of fixed discriminant $\Delta$.

\noindent Let $B(N,p)=\{-\Delta N,p\}=\QQ+\QQ i^N+\QQ j^N+\QQ k^N$ be a representation of $B$; we will write $R(N)\subset B(N,p)$ to denote the Eichler order of level $N$ of Hashimoto \cite{Hashimoto95}: $R(N)=\ZZ e_1^N+\ZZ e_2^N+\ZZ e_3^N+\ZZ e_4^N$ with  $$e_1^N=1,\  e_2^N=\frac{1+j^N}{2},\ e_3^N=\frac{i^N+k^N}{2},\ e_4^N=\frac{a\Delta Nj^N+k^N}{p}$$ where $a\in\ZZ$ satisfies $a^2\Delta N+1\equiv\ 0\ \modulo\ p$.  By abuse of notation, in this section we will write $R(M)$ instead of $\Psi_M^N(R(M))$.  In this way, if $N|M$ the inclusion $R(M)\subset R(N)$ in $B(N,p)$ is true.

\begin{lemma}\label{rango}
Let $B$ be a quaternion algebra over $\QQ$ of discriminant $\Delta$; let $N$ be a positive integer prime to $\Delta$ and $q$ be a prime number not dividing $\Delta$. Then the $\ZZ$-rank of $\bigcap_{n\in\NN}R(Nq^n)$ is equal to the $\ZZ$-rank of $\bigcap_{n\in\NN}R(q^n)$.
\end{lemma}
\proof
Let $B(1,p)$ be a representation of $B$ where $p$ is as in Section \ref{sect:preliminari}. It is obvious that 
\begin{equation}\label{immaginezeroN}
\bigcap_nR(Nq^n)=\bigcap_nR(q^n)\cap R(N)\subset R(1).
\end{equation}
Since the rank is invariant by isomorphism and $R(N)$ has maximal rank over $\ZZ$, then $$\textnormal{rk}\left(\bigcap_{n\in\NN}R(Nq^n)\right)=\textnormal{rk}\left(\bigcap_{n\in\NN}R(q^n)\right).$$\qed

\noindent Let $B(1,p)$ be a representation of $B$ anf let $q\not|\Delta$ be a prime number; we consider the chain of Eichler orders $$...\subset R(q^n)\subset...\subset R(q^2)\subset R(q)\subset R(1)$$ in $B(1,p)$. We will characterize the intersection $\mathcal A_q=\bigcap_{n\in\NN}R(q^n)$ as $\ZZ$-lattice. Since $$R(q)\simeq R(1)\cap(\varphi_q^1)^{-1}\left(
\begin{array}{ll}
\ZZ_q & \ZZ_q \\
q\ZZ_q & \ZZ_q 
\end{array}
\right)$$  where $\varphi_q^1:B(1,p)_q\to M_2(\QQ_q)$ is a local isomorphism, then 
\begin{equation}\label{immaginezero}
\mathcal A_q
\simeq R(1)\cap\left[\bigcap_n(\varphi_q^1)^{-1}\left(
\begin{array}{ll}
\ZZ_q & \ZZ_q \\
q^n\ZZ_q & \ZZ_q 
\end{array}
\right)\right]=R(1)\cap(\varphi_q^1)^{-1}\left(
\begin{array}{ll}
\ZZ_q & \ZZ_q \\
0 & \ZZ_q 
\end{array}
\right).
\end{equation}

\begin{proposition}
Let $B$ be a quaternion algebra over $\QQ$ of discriminant $\Delta$; let we fix a prime number $q$ not dividing $\Delta$. The intersection $$\mathcal A_q=\bigcap_{n\in\NN} R(q^n)$$ has rank 2 over $\ZZ$.
\end{proposition}
\proof\\
Let we fix a prime number $p$ is as in Section \ref{sect:preliminari}. 
We will distinguish the following cases:
\begin{enumerate}
\item $q\not|\Delta p$ such that $\left(\frac{p}{q}\right)=-1$;
\item $q\not|\Delta$ such that $\left(\frac{p}{q}\right)=1$;
\item $q=p$.
\end{enumerate}

\noindent 1. Let $q\not|\Delta$ be a prime number such that $\left(\frac{p}{q}\right)=-1$. Let $N$ be a positive integer prime to $\Delta$ such that $\left(\frac{p}{s}\right)=1$ for all $s|N$ and $\left(\frac{-\Delta N}{q}\right)=1$. This last condition on $N$ implies that there exists $x(N,q)\in\ZZ_q$ such that $-\Delta N=x(N,q)^2$. Let we represent $B$ as $B(N,p)=\{-\Delta N,p\}$.\\
If $q$ is such that $-\Delta$ is a square in $\ZZ_q$, then we can take $N=1$ and $\mathcal A_q\subset R(1)$ in $B(1,p)$;  if $h\in\mathcal A_q$ then by (\ref{immaginezero}) there exist $\alpha,\beta,\gamma,\delta\in\ZZ$ such that $h=\alpha e_1^1+\beta e_2^1+\gamma e_3^1+\delta e_4^1$ where $\{e_1^1,e_2^1,e_3^1,e_4^1\}$ is the Hashimoto basis of $R(1)$ in $B(1,p)$. Moreover, by the identity (\ref{q^2|N}) $$x(N,q)\left[-\frac{\gamma}{2}-\frac{\delta}{p}\right]+\frac{\beta}{2}+\frac{\delta a\Delta}{p}=0.$$ Then  $\mathcal A_q\subset R(1)$ can be expressed as the $\ZZ$-lattice $\mathcal A_q=\ZZ e_1^1+\ZZ e^1$ where $e^1=-2a\Delta e_2^1-2e_3^1+pe_4^1$.\\
If $q$ is such that $-\Delta$ is not a square in $\ZZ_q$, then we take $N$ such that $\left(\frac{N}{q}\right)=-1$. If $h\in\bigcap_nR(Nq^n)\simeq R(N)\cap(\varphi_q^N)^{-1}\left(
\begin{array}{ll}
\ZZ_q & \ZZ_q \\
0 & \ZZ_q 
\end{array}
\right)$ in $B(N,p)$ then
there exist $\alpha,\beta,\gamma,\delta\in\ZZ$ such that $h=\alpha e_1^N+\beta e_2^N+\gamma e_3^N+\delta e_4^N$ where $\{e_1^N,e_2^N,e_3^N,e_4^N\}$ is the Hashimoto basis of $R(N)$ in $B(N,p)$. Moreover, by the identity (\ref{q^2|N}) $$x(N,q)\left[-\frac{\gamma}{2}-\frac{\delta}{p}\right]+\frac{\beta}{2}+\frac{\delta a\Delta N}{p}=0.$$ Then $\bigcap_nR(Nq^n)\subset R(N)$ can be expressed as the $\ZZ$-lattice $\bigcap_nR(Nq^n)=\ZZ e_1^N+\ZZ e^N$ where $e^N=-2a\Delta N e_2^1-2e_3^1+pe_4^1$. By Lemma \ref{rango}, $\textnormal{rk}\left(\bigcap_nR(q^n)\right)=\textnormal{rk}\left(\bigcap_nR(Nq^n)\right)=2$.

\noindent 2. Let $q\not|\Delta$ be a prime number such that $\left(\frac{p}{q}\right)=1$; we represent the quaternion algebra $B$ as $B(1,p)=\{-\Delta,p\}$. If $h\in\mathcal A_q$, then by (\ref{immaginezero})  and by the identity (\ref{eldue}), $h=\alpha e_1^1+\beta e_2^1+\gamma e_3^1+\delta e_4^1$ where 
$\alpha,\beta,\gamma,\delta\in\ZZ$ satisfy the  equation 
\begin{equation}\label{condzero2}
\sqrt{p}(-\gamma)+(\gamma p+2\delta)=0.
\end{equation} 
Then $\gamma=\delta=0$ and $\mathcal A_q\subset R(1)$ can be expressed as the $\ZZ$-lattice $\mathcal A_q=\ZZ e_1^1+\ZZ e_2^1$.

\noindent 3. Let $q=p$; we represent the quaternion algebra $B$ as $B(1,p)=\{-\Delta,p\}$. If $h\in\mathcal A_p$, then by (\ref{immaginezero}) and by the identity (\ref{p=q}), $h=\alpha e_1^1+\beta e_2^1+\gamma e_3^1+\delta e_4^1$
where $\alpha,\beta,\gamma,\delta\in\ZZ$ satisfy the  equation 
\begin{equation}\label{condzero3}
\sqrt{-\Delta}(\gamma p+2\delta)+(\beta p+2a\Delta\delta)=0.
\end{equation}
Then $\mathcal A_q\subset R(1)$ can be expressed as the $\ZZ$-lattice $\mathcal A_p=\ZZ e_1^1+\ZZ e$ where $e=-2a\Delta e_2^1-2e_3^1+pe_4^1$.
\qed
{}

\begin{proposition}
Let $B$ a quaternion algebra over $\QQ$ of discriminant $\Delta$ and let $B(1,p)$ be a representation of $B$. Let $q,s\not|\Delta p$ be two prime number  such that $\left(\frac{p}{q}\right)=1$ and $\left(\frac{p}{s}\right)=-1$. Then $$\mA_q\cap\mA_p=\mA_s\cap\mA_p=\mA_q\cap\mA_s=\ZZ.$$
\end{proposition}
\proof
Let $\{e_1, e_2, e_3, e_4\}$ be the Hashimoto basis  of $R(1)$ in $B(1,p)$.\\
If $h\in\mA_q\cap\mA_p$ then $h=\alpha e_1+\beta e_2+\gamma e_3+\delta e_4$ where  $\alpha,\beta,\gamma,\delta\in\ZZ$ satisfy the equations (\ref{condzero2}) and (\ref{condzero3}). This imply that $\beta=\gamma=\delta=0.$\\
If $h\in\mA_s\cap\mA_p$ then $h=\alpha e_1+\beta e_2+\gamma e_3+\delta e_4$ where $\alpha,\beta,\gamma,\delta\in\ZZ$ satisfy the equation  (\ref{condzero3}) and by the identity (\ref{q^2|N}) 
\begin{equation}\label{zero0}
x\left[-\frac{\gamma}{2}-\frac{\delta}{p}\right]+y\left[\frac{\gamma}{2}\right]+\frac{\beta}{2}+\delta\frac{a\Delta}{p}=0.
\end{equation} where $x,y\in\ZZ_q$ are such that $-\Delta=x^2-py^2$. This imply that $\beta=\gamma=\delta=0.$\\
If $h\in\mA_q\cap\mA_s$ then $h=\alpha e_1+\beta e_2+\gamma e_3+\delta e_4$ where  $\alpha,\beta,\gamma,\delta\in\ZZ$ satisfy the equations (\ref{condzero2}) and (\ref{zero0}). This imply that $\beta=\gamma=\delta=0.$ \qed

\begin{corollary}
Let $B$ a quaternion algebra over $\QQ$ of discriminant $\Delta$; then $$\mA:=\bigcap_{N}R(N)=\ZZ$$ where $N$ runs over the set of positive integer numbers primes to $\Delta$. 
\end{corollary}

\noindent As corollary, by Lemma \ref{uno}, the following result holds:  
\begin{corollary}
Let $\Phi(N)$ be the group  defined in Section \ref{caratterizzazione}, then:
 $$\bigcap_{N}\Phi(N)=\{\pm 1\}$$ where $N$ runs over the set of positive integer numbers primes to $\Delta$.
\end{corollary}

\section{Example}

Using a mathematical problem-solving environment as Maple, wich work with $p$-adic numbers, it is possible to produce some examples.\\


\noindent Let we consider the quaternion algebra $B$ over $\QQ$ with discriminant $\Delta=35$; following Hashimoto we can represent it as $B(3,13)=\{-105,13\}$.
A basis over $\ZZ$ of the Eichler order $R(3)$ of $B(3,13)$ is $$e_1=1,\ e_2=\frac{1+j}{2},\ e_3=\frac{i+k}{2},\ e_4=\frac{525j+k}{13}.$$ 
If we consider $q=11$, then $q\not|\Delta$ and $p=13$ is not a square in $\ZZ_{11}^\s$; thus by Proposition \ref{NNqualunque}, a basis of the Eichler order $R(33)$ in $B(3,13)$ is: 
$$f_1=1,\ \ \  f_2=-\frac{5}{2}+\frac{i}{2}-\frac{5}{2}j+\frac{k}{2}$$ $$f_3=-3150-\frac{40425}{13}j+\frac{1}{13}k,\ \ \ f_4=\frac{11}{2}+\frac{11}{2}j.$$ A basis of  $\delta_{11}R(33)\delta_{11}^{-1}$ in $B(3,13)$ is: $$g_1=1,\ \ \ g_2=\frac{5}{2}+\frac{i}{2}+\frac{5}{2}j+\frac{k}{2}$$ $$g_3=-3150-\frac{40425}{13}j+\frac{1}{13}k\ \ \ g_4=\frac{11}{2}+\frac{11}{2}j.$$

\noindent Moreover let $B(17,13)=\{-595,13\}=\QQ+\QQ i^{(17)}+\QQ j^{(17)}+\QQ k^{(17)}$ be the quaternion algebra over $\QQ$ with discriminant $35$ and $N=17$; then $(i^{(17)})^2=-595$ and $(j^{(17)})^2=13$. It is possible to write the isomorphism $\Psi^{17}_3:B(3,13)\to B(17,13)$ described in Section \ref{isomo}:
$$\Psi^{17}_3(j)=j^{(17)},\ \ \ \ \ \Psi^{17}_3(i)=\frac{8}{17}i^{(17)}+\frac{1}{17}k^{(17)}.$$

\bigskip

\begin{center}
\noindent {\bf Authors' affiliation:}\\

\bigskip

\noindent Miriam Ciavarella\\ Universit\`a degli Studi di Torino\\ Dipartimento di Matematica\\ Via Carlo Alberto,10\\ 10123 Torino (Italy)\\ e-mail: miriam.ciavarella@unito.it\\

\bigskip

\noindent Lea Terracini\\ Universit\`a degli Studi di Torino\\ Dipartimento di Matematica\\ Via Carlo Alberto,10\\ 10123 Torino (Italy)\\ e-mail: lea.terracini@unito.it
\end{center}

\end{document}